\documentclass[12pt]{article}

\usepackage[left=2.5cm, right=2.5cm, top=3cm, bottom=3cm]{geometry}

\usepackage{hyperref}
\usepackage{url}
\usepackage{algorithm}
\usepackage{algorithmic}
\usepackage{amsthm,amsmath,amssymb}
\usepackage{amsfonts,dsfont} 
\usepackage{mathtools} 
\usepackage{graphics}
\usepackage{epsfig}
\usepackage{tablefootnote}
\usepackage{cite}
\usepackage{authblk}
\newtheorem{theo}{Theorem}[section]  

\newtheorem{lemm}[theo]{Lemma}

\newcommand{\beq}{\begin{equation}}
\newcommand{\eeq}{\end{equation}}
\newcommand{\beqa}{\begin{eqnarray}}
\newcommand{\eeqa}{\end{eqnarray}}
\newcommand{\beqas}{\begin{eqnarray*}}
\newcommand{\eeqas}{\end{eqnarray*}}
\newcommand{\beqs}{\begin{equation*}}
\newcommand{\eeqs}{\end{equation*}}

\newcommand{\prox}{\text{Prox}}

\title{On the Convergence of Primal-Dual Proximal Incremental Aggregated Gradient Algorithms}
\author{Xianchen Zhou \thanks{NUDT. email: zhouxianchen13@nudt.edu.cn} ~	
Wei Peng \thanks{NUDT. email: weipeng0098@126.com.} ~Hongxia Wang \thanks{NUDT. email: wanghongxia@nudt.edu.cn}}

\date{\today}
\begin{document}

\maketitle

\begin{abstract} 
In this paper, we adapt proximal incremental aggregated gradient methods to saddle point problems, which is motivated by decoupling linear transformations in regularized empirical risk minimization models. First, the Primal-Dual Proximal Incremental Aggregated (PD-PIAG) methods with extrapolations were proposed. We proved that the primal-dual gap of the averaged iteration sequence sublinearly converges to 0, and the iteration sequence converges to some saddle point. Under the strong convexity of $f$ and $h^\ast$, we proved that the iteration sequence linearly converges to the saddle point. Then, we propose a PD-PIAG method without extrapolations. The primal-dual gap of the iteration sequence is proved to be sublinearly convergent under strong convexity of $f$. 
\end{abstract}

\section{Introduction}
We consider the convex-concave saddle point problems of the form
\begin{align}
\min_{x\in\mathbb{R}^{d_1}}\max_{y\in\mathbb{R}^{d_2}}\mathcal{L}(x,y):=f(x)+\langle Kx,y\rangle-h^\ast(y),\label{pm2}
\end{align}
where $f:=\sum_{i=1}^M f_i$ with $f_i:\mathbb{R}^{d_1}\rightarrow \mathbb{R}$ being smooth and convex, conjugate function $h^\ast:\mathbb{R}^{d_2}\rightarrow \mathbb{R}$ is convex and possibly nonsmooth, and $K\in \mathbb{R}^{d_2\times d_1}$ is a matrix. We are motivated by the optimization problems with a cost function that consists of additive components and a regularizer:
\begin{align}
\min_{x}\sum_{i=1}^M f_i(x)+h(Kx),\label{pm1}
\end{align}
Consider the conjugate function
\begin{align}
h(Kx)=\max_{y\in\mathbb{R}^{d_2}} \left\{\langle\label{pm44} K^Ty,x\rangle-h^\ast(y)\right\},
\end{align}
Substituting \eqref{pm44} into \eqref{pm1}, we obtain the fundamental optimization model \eqref{pm2}, which emerges in numerous problems, including machine learning, signal processing, imaging science, communication systems, and distributed optimization. 

For the cases where $K=I$ and the proximal operation of $h$(or $h^\ast$) is inexpensive, \eqref{pm1} can be handle by the well-known forward-backward splitting (FBS) \cite{combettes2005signal}:
\begin{align}
x_{k+1}=\arg\min_x\left\{ h(x)+\left\langle  \sum_{i=1}^M \nabla f_i(x_k),x-x_k\right\rangle+\frac{\tau_k}{2}\|x-x_k\|^2\right\}.
\end{align}
In many cases, the number of component functions $M$ is so large that directly computing the full gradient of the smooth part becomes prohibitive. To overcome this difficulty, many stochastic variants of FBS have been proposed\cite{xiao2014proximal,nitanda2014stochastic,defazio2014saga}. Besides these stochastic methods, the proximal incremental aggregated gradient (PIAG) method, as a deterministic method, is presented:
\begin{align}
x_{k+1}=\arg\min_x\left\{ h(x)+\left\langle  \sum_{i=1}^M \nabla f_i(x_{k-\tau_k^i}),x-x_k\right\rangle+\frac{\tau_k}{2}\|x-x_k\|^2\right\},
\end{align}
where $\tau_k^i$ is the delay of the $k$-th iteration of the $i$-th component.

The key idea of PIAG is to construct an ``inexact gradient" to substitute for the full gradient at each iteration. 
PIAG has been extensively investigated recently under the strong convex assumption, and its global linear convergence of the objective function and iterative sequences have been established \cite{RN63,PIAG2016}.
A PIAG-like algorithmic framework has been proposed, which includes PIAG as a special  case \cite{zh1}, and a linear convergence theory was built, but under strictly weaker assumptions. 
The linear convergence of nonconvex PIAG under error bound conditions has also been studied in \cite{peng}.

However, all current works of PIAG are aimed to minimize the special cases where $K=I$.
Though the proximal operation for $h$ can be inexpensive, its composition with a linear transformation $K$ may be too expensive. Therefore, we wonder if  incremental aggregated methods could be applied to \eqref{pm2}.

There exists a large amount of literature on primal-dual algorithms to solve \eqref{pm2}. Gradient methods for solving saddle point problems have attracted much research since the seminal work of \cite{arrow1958studies}. \cite{nedic2009subgradient} extended their method, and performs subgradient steps on the primal and dual variables alternatingly. Many variants were proposed \cite{chambolle2016ergodic,esser2010general,chernov2016fast}. Typically, under strongly convexity of $f$ and $h^\ast$, linear convergence can be guaranteed \cite{chambolle2011first}. Preconditioned and adaptive stepsizes versions were proposed in \cite{pock2011diagonal,combettes2014variable,goldstein2013adaptive}.  We refer readers to \cite{komodakis2015playing} for a detailed review.

\textbf{Contributions. }Our main contribution is to adapt PIAG to Primal-Dual methods and obtain PD-PIAG, which can be implemented by an asynchronous distributed framework. First we study PD-PIAG with extrapolations, linear convergence and sublinear convergence of which are obtained under different assumptions respectively. Then we study an Arrow-Hurwicz method like version of incremental aggregated gradient methods with a sublinear convergence guarantee.
\section{Notations, Preliminaries and Algorithms}
Throughout the paper, $d$-dimensional Euclidean space is denoted by $\mathbb{R}^d$ and its inner product by $\langle\cdot,\cdot \rangle$. The $l_2$-norm is denoted by $\|\cdot\|$. 
The gradient operator of a differentiable function is denoted by $\nabla$.
The subdifferential of a proper closed convex function $G$ is defined by
\begin{align}
\partial G(x):=\left\{v\in\mathbb{R}^d: G(u)-G(x)-\langle v,u-x\rangle\geq0,\forall u\in\mathbb{R}^d\right\}
\end{align}
The proximal operator of a proper closed function $G$ is defined by
\begin{align}
\text{Prox}_G(\cdot):=\arg\min_{x\in\mathbb{R}^d}\left\{G(x)+\frac{1}{2}\|x-\cdot\|^2\right\}
\end{align}

\subsection{Saddle points and Min-Max problem}
We consider the saddle-point problem:
\begin{align}
\min_{x\in \mathbb{R}^{d_1}}\max_{y\in\mathbb{R}^{d_2}} \mathcal{L}(x,y).\label{saddle}
\end{align}
We say $(\hat x, \hat y)$ is a saddle point for $\mathcal{L}$ if
\begin{align}
\mathcal{L}(\hat x,y)\leq \mathcal{L}(\hat x,\hat y)\leq \mathcal{L}(x,\hat y),~~~\forall(x,y) \in\mathbb{R}^{d_1}\times \mathbb{R}^{d_2}.\label{equiv1}
\end{align}
If $(\hat x,\hat y)$ is a saddle point, then we have
\begin{align}
\mathcal{L}(x,\hat y)-\mathcal{L}(\hat x,y)\geq 0, ~~~\forall (x,y)\in\mathbb{R}^{d_1}\times\mathbb{R}^{d_2}.\label{lieq}
\end{align}
If $f$ and $g^\ast$ are convex, proper, and closed, then \eqref{equiv1} is equivalent to
\begin{align}
0\in K^T\hat y+\partial f(\hat x),~~~0\in K\hat x-\partial h^\ast(\hat y). \label{opt1}
\end{align}
We introduce the partial primal-dual gap\cite{chambolle2011first}. For closed set $B_1\subset\mathbb{R}^{d_1}$ and $B_2\subset\mathbb{R}^{d_2}$, define
\begin{align}
\mathcal{G}_{B_1\times B_2}(x,y):=&\max_{y'\in B_2}\mathcal{L}(x,y')-\min_{x'\in B_1}\mathcal{L}(x',y).
\end{align}
If $(x,y)\in B_1\times B_2$, we have
\begin{align}
\mathcal{G}_{B_1\times B_2}(x,y)=&\max_{y'\in B_2}[\mathcal{L}(x,y')-\mathcal{L}(x,y)]-\min_{x'\in B_1}[\mathcal{L}(x',y)-\mathcal{L}(x,y)]\nonumber\\
\geq &[\mathcal{L}(x,y)-\mathcal{L}(x,y)]-[\mathcal{L}(x,y)-\mathcal{L}(x,y)]=0.
\end{align}
Conversely, if $\mathcal{G}_{B_1\times B_2}(x,y)=0$, we have
\begin{align}
y\in \arg\max_{y'\in B_2}[\mathcal{L}(x,y')-\mathcal{L}(x,y)],~~~x\in \arg\min_{x'\in B_1}[\mathcal{L}(x',y)-\mathcal{L}(x,y)].\label{opt2}
\end{align}
Furthermore, if $(x,y)$ is an interior point of $B_1\times B_2$, then \eqref{opt2} leads to
\begin{align*}
0\in K^Ty+\partial f(x),~~~0\in Kx-\partial h^\ast(y),
\end{align*}
which implies $(x,y)$ is a saddle point according to \eqref{opt1}. Therefore, if $B_1\times B_2$ is large enough, then $\mathcal{G}_{B_1\times B_2}(x,y)$ can measure the optimality of $(x,y)$.

\subsection{Assumptions}
We listed some assumptions to be used in this manuscript as follows:
\begin{enumerate}
	\item [A1] For $1\leq i \leq N$, $f_i$ is $L_i$-smooth, i.e.,
	\begin{align}
	\left|f_i(y)-f_i(x)-\langle \nabla f_i(x),y-x\rangle \right|\leq L_i\frac{\|y-x\|^2}{2},
	\end{align}
	where $L := \sum_{i=1}^N L_i>0$.
	
	\item [A2] $h^\ast:\mathbb{R}^{d_2}\rightarrow(-\infty,\infty]$ is proper, closed, and convex.
	
	\item [A3] The time-varing delays $\tau_k^i$ are bounded, i.e., there exists a nonnegative integer $T$ such that $\forall k\geq 1,i\in\{1,2,\cdots,M\}$, we have
	\begin{align}
	\tau_{k}^i\in\{0,1,\cdots,T\},
	\end{align}
	where $T$ is called the \textit{delay parameter}.
	The above three assumptions are some standard assumptions in PIAG research, Two assumptions  are established to get better convergence.
	\item [B1] For $1\leq i \leq N$, $f_i$ satisfies
	\begin{align}
	f_i(y)\geq f_i(x)+\langle \nabla f_i(x),y-x\rangle+\delta_i\frac{\|y-x\|^2}{2},
	\end{align}
	Define $\delta := \sum_{i=1}^M \delta_i>0$, where $\delta_i$ is not required to be nonnegative.
	\item [B2] $h^\ast$ is $\gamma$-strongly convex, i.e.,
	\begin{align}
	h^\ast(x)\geq g^\ast(x)+\langle v,y-x\rangle+\gamma \frac{\|y-x\|^2}{2},~~~\forall v\in \partial h^\ast(x),
	\end{align}
	where $\gamma>0$.
\end{enumerate}
Assumption B2 is equivalent to $1/\gamma$ smoothness of $h$.

\subsection{Algorithms}
Similar to PIAG, we consider PD-PIAG fomation:
\begin{equation} \label{algorithm1}
\text {(PD-PIAG)}\left\{\begin{aligned} g_{k} &=\sum_{i=1}^{N} \nabla f_{i}\left(x_{k-\tau_{k}^{i}}\right) \\ 
x_{k+1} &=x_{k}-\sigma\cdot g_{k}-\sigma\cdot K^T\bar y \\
y_{k+1} &=\prox_{\tau h^\ast}\left(y_k+\tau\cdot Kx_{k+1}\right), \end{aligned}\right.
\end{equation}
First line is the aggregation of $ \nabla f_{i}$, second line is the gradient update of $x$ and second line is the gradient update of proximity operator. $\bar y$ is an undermined variable. Different choice of  $\bar y$ can lead to different convergence, and we discuss convergence in next section.

Here algorithm \ref{Alg:pd_piag_cycle} proposes  an implementation of Update gradient circularly PD-PIAG.

\begin{algorithm}
	\caption{Update gradient circularly~PD-PIAG~\label{Alg:pd_piag_cycle}}
	\begin{algorithmic}
		\STATE{\textbf{Require:} 初始点~$(x_0,y_0)\in\mathbb{R}^{d_1}\times\mathbb{R}^{d_2},\{e_i:=\nabla f_i(x_0)\}_{i=1}^N$}
		\STATE{\textbf{Initialization:}~$g_0:=\sum_{i=1}^N e_i,k:=0$}
		\STATE{1: \textbf{repeat}}
		\STATE{2: \qquad ~$i_k:=(k \mod N)+1$}
		\STATE{3: \qquad~$x_{k+1}:=x_k-\sigma g_k-\sigma K^T \bar y$}
		\STATE{3: \qquad ~$y_{k+1}:=\prox_{\tau h^\ast}(y_k+\tau Kx_{k+1})$}
		\STATE{5: \qquad Compute~$\nabla f_{i_k} (x_{k+1})$}
		\STATE{6: \qquad Update sum of gradients~$g_{k+1}:=g_k+\nabla f_{i_k}(x_{k+1})-e_{i_k}$}
		\STATE{7: \qquad Update memory~$e_{i_k}:=\nabla f_{i_k} (x_{k+1})$}
		\STATE{8: \qquad $k:=k+1$}
		\STATE{9: \textbf{until:} termination condition satisfied}
		\STATE{\textbf{Return:} $\{(x_{k}, y_{k})\}$}
	\end{algorithmic}
\end{algorithm}

\section{Convergence Analysis}

\subsection{Sublinear convergence}
In the first analysis of convergence , take the  $\bar y = 2y_{k}-y_{k-1}$ in PD-PIAG. When iteration stepsize $\sigma$ and $\gamma$ are small enough, sublinear convergence of partial primal-dual gap $\mathcal{G}_{B_1\times B_2}(x,y)$ can obtained.
\begin{theo}\label{theo1} Assume A1-A3 hold. Take $\bar y = 2y_{k}-y_{k-1}$  in each iteration. Assume the problem has a saddle point $(\hat x,\hat y)$. Choose $\tau$ and $\sigma$ such that 
\begin{align}
\sqrt{\tau\sigma}\|K\|+\sigma L(T+1)^2<1.\label{cond1}
\end{align}
Then:
\begin{enumerate}
\item[\rm(i)]
 The sequence $\{x_k\},\{y_k\}$ is bounded since
\begin{align}
\frac{\left\|x_{n}-\hat{x}\right\|^{2}}{2 \sigma}+\frac{\left\|y_{k}-\hat{y}\right\|^{2}}{2 \tau}
\quad \leq C\left(\frac{\left\|x_{0}-\hat{x}\right\|^{2}}{2 \sigma}+\frac{\left\|y_{0}-\hat{y}\right\|^{2}}{2 \tau}\right),
\end{align}
where the constant $C=(1-\tau\sigma\|K\|^2)^{-1}$. 
\item[\rm(ii)]
Define the averaged sequences $\bar x_M=(\sum_{k=1}^M x_k)/M$ and $\bar y_M=(\sum_{k=1}^M y_k)/M$ for all $M>0$.
Then for any bounded closed set $B_1\times B_2\subset X\times Y$, the restricted gap has the following bound
\begin{align}
\mathcal{G}_{B_{1} \times B_{2}}\left(\bar x_{M}, \bar y_{M}\right) \leq \frac{1}{M}\max _{(x, y) \in B_{1} \times B_{2}} \left\{\frac{\left\|x-x_{0}\right\|^{2}}{2 \sigma}+\frac{\left\|y-y_{0}\right\|^{2}}{2 \tau}\right\}.
\end{align}
Moreover, the cluster points of $\{(x_M,y_M)\}$ are saddle points. 
\item[\rm(iii)]
There exists a saddle point $(x^\ast,y^\ast)$ such that $x_k\rightarrow x^\ast$ and $y_	k\rightarrow y^\ast$.
\end{enumerate}
\end{theo}

\textbf{Proof.}
We first introduce the descent property of $f$. Via the $L_i$-smoothness of $f_i$, we have
\begin{align}
f_i (x)&\geq f_i({x_{k-\tau_k^i}})+\left\langle \nabla f_i(x_{k-\tau_k^i}), x-x_{k-\tau_k^i} \right\rangle\nonumber\\
&\geq f_i({x_{k+1}})+\left\langle \nabla f_i (x_{k-\tau_k^i}), x-x_{k+1} \right\rangle-\frac{L_i\|x_{k-\tau_k^i}-x_{k+1}\|^2}{2}\nonumber\\
&\geq f_i({x_{k+1}})+\left\langle \nabla f_i(x_{k-\tau_k^i}), x-x_{k+1} \right\rangle-\frac{L_i(T+1)}{2}\sum_{j=k-T}^k\|x_{j+1}-x_{j}\|^2\label{f2}
\end{align}
Summing \eqref{f2} from $i=1$ to $N$ yields
\begin{align}
f(x) \geq & f\left(x_{k+1}\right)+\left\langle g_k, x-x_{k+1}\right\rangle-\frac{L(T+1)}{2}\sum_{j=k-T}^k\|x_{j+1}-x_j\|^2\nonumber \\
=&f\left(x_{k+1}\right)+\left\langle \frac{x_{k}-x_{k+1}}{\sigma}-K^T \overline{y}, x-x_{k+1}\right\rangle-\frac{L(T+1)}{2}\sum_{j=k-T}^k\|x_{j+1}-x_j\|^2.\label{f}
\end{align}
According to the iteration procedure, we have
\begin{align}
\partial h^\ast\left(y_{k+1}\right) \ni& \frac{y_{k}-y_{k+1}}{\tau}+K x_{k+1}\label{g3}.
\end{align}
Combining convexity of $g^\ast$ and \eqref{g3}, we obtain
\begin{align}  
h^\ast(y) \geq & h^\ast\left(y_{k+1}\right)+\left\langle\frac{y_{k}-y_{k+1}}{\tau}, y-y_{k+1}\right\rangle +\left\langle y-y_{k+1}, Kx_{k+1}\right\rangle. \label{g}\end{align}
Summing \eqref{f} and \eqref{g}, it follows that
\begin{align}
\frac{\left\|x-x_{k}\right\|^2}{2 \sigma}+\frac{\left\|y-y_{k}\right\|^{2}}{2 \tau}
\geq&\mathcal{L}(x_{k+1},y)-\mathcal{L}(x,y_{k+1})+\frac{\left\|x-x_{k+1}\right\|^{2}}{2 \sigma}+\frac{\left\|y-y_{k+1}\right\|^{2}}{2 \tau}\nonumber\\
&+\frac{\left\|x_{k}-x_{k+1}\right\|^{2}}{2 \sigma}+\frac{\left\|y_{k}-y_{k+1}\right\|^{2}}{2 \tau}-\frac{L(T+1)}{2}\sum_{j=k-T}^k\|x_{j+1}-x_j\|^2\nonumber\\
&-\left\langle K^T\left(y_{k+1}-\overline{y}\right), x_{k+1}-x\right\rangle.\label{ab}
\end{align}
Substituting $\bar y=2y_{k}-y_{k-1}$ for the last term of \eqref{ab}, we have
\begin{align*}
&-\left\langle K^T\left(y_{k+1}-\overline{y}\right), x_{k+1}-x\right\rangle\\
=&-\left\langle K^T\left(\left(y_{k+1}-y_{k}\right)-\left(y_{k}-y_{k-1}\right)\right), x_{k+1}-x\right\rangle\\
=&-\left\langle K^T\left(y_{k+1}-y_{k}\right), x_{k+1}-x\right\rangle+\left\langle K^T\left(y_{k}-y_{k-1}\right), x_{k}-x\right\rangle +\left\langle K^T\left(y_{k}-y_{k-1}\right), x_{k+1}-x_{k}\right\rangle \\
\geq&-\left\langle K^T\left(y_{k+1}-y_{k}\right), x_{k+1}-x\right\rangle+\left\langle K^T\left(y_{k}-y_{k-1}\right), x_{k}-x\right\rangle - \|K\|\left\|y_{k}-y_{k-1}\right\|\left\|x_{k+1}-x_{k}\right\|.\end{align*}
Since fundamental inequality
\begin{align} \|K\|\left\|x_{k+1}-x_{k}\right\|\left\|y_{k}-y_{k-1}\right\| \leq & \sqrt{\tau\sigma}\|K\|\frac{\left\|x_{k+1}-x_{k}\right\|^{2}}{2\sigma} +\sqrt{\tau\sigma}\|K\|\frac{\left\|y_{k}-y_{k-1}\right\|^{2}}{2\tau}, \end{align}
we have
\begin{align}
\frac{\left\|x-x_{k}\right\|^2}{2 \sigma}+\frac{\left\|y-y_{k}\right\|^{2}}{2 \tau}
\geq&\mathcal{L}(x_{k+1},y)-\mathcal{L}(x,y_{k+1})\nonumber\\
&+\frac{\left\|x-x_{k+1}\right\|^{2}}{2 \sigma}+\frac{\left\|y-y_{k+1}\right\|^{2}}{2 \tau}+(1-\sqrt{\tau\sigma}\|K\|)\frac{\left\|x_{k}-x_{k+1}\right\|^{2}}{2 \sigma}\nonumber\\
&+\frac{\left\|y_{k}-y_{k+1}\right\|^{2}}{2 \tau}-\sqrt{\tau\sigma}\|K\|\frac{\left\|y_{k-1}-y_{k}\right\|^{2}}{2 \tau}\nonumber\\
&-\left\langle K^T\left(y_{k+1}-y_k\right), x_{k+1}-x\right\rangle+\left\langle K^T\left(y_{k}-y_{k-1}\right), x_{k}-x\right\rangle\nonumber\\
&-\frac{L(T+1)}{2}\sum_{j=k-T}^k\|x_{j+1}-x_j\|^2\label{fs491}.
\end{align}
Summing \eqref{fs491} from $k=0$ to $M-1$, we obtain
\begin{align}
\frac{\left\|x-x_{0}\right\|^2}{2 \sigma}+\frac{\left\|y-y_{0}\right\|^{2}}{2 \tau}
\geq&\sum_{k=1}^{M}\left(\mathcal{L}(x_k,y)-\mathcal{L}(x,y_k)\right)\nonumber\\
&+\frac{\left\|x-x_{M}\right\|^{2}}{2 \sigma}+\frac{\left\|y-y_{M}\right\|^{2}}{2 \tau}+(1-\sqrt{\tau\sigma}{\|K\|})\sum_{k=0}^{M-1}\frac{\left\|x_{k}-x_{k+1}\right\|^{2}}{2 \sigma}\nonumber\\
&+\left(1-{\sqrt{\tau\sigma}\|K\|}\right)\sum_{k=1}^{M-1}\frac{\left\|y_{k-1}-y_{k}\right\|^{2}}{2 \tau}+\frac{\left\|y_{M}-y_{M-1}\right\|^{2}}{2 \tau}\nonumber\\
&-\frac{L(T+1)^2}{2}\sum_{k=0}^{M-1}\|x_{k+1}-x_k\|^2-\left\langle K^T\left(y_{M}-y_{M-1}\right), x_{M}-x\right\rangle.\label{ref73}
\end{align}
As before, we have
\begin{align}
\left | \left\langle K^T(y_{M}-y_{M-1}),x_{M}-x)\right\rangle \right |&\leq \|K\|\cdot\|y_{M}-y_{M-1}\|\cdot\|x_{M}-x\|\nonumber\\
&\leq\frac{\|y_{M}-y_{M-1}\|^2}{2\tau}+\tau\sigma\|K\|^2\frac{\|x_{M}-x\|^2}{2\sigma}.\label{ref74}
\end{align}
Combining \eqref{ref73} and \eqref{ref74}, we have
\begin{align}
\frac{\left\|x-x_{0}\right\|^2}{2 \sigma}+\frac{\left\|y-y_{0}\right\|^{2}}{2 \tau}
\geq&\sum_{k=1}^M\left(\mathcal{L}(x_k,y)-\mathcal{L}(x,y_k)\right)\nonumber\\
&+\frac{\left\|x-x_{M}\right\|^{2}}{2 \sigma}+\frac{\left\|y-y_{M}\right\|^{2}}{2 \tau}+(1-\sqrt{\tau\sigma}\|K\|)\sum_{k=0}^{M-1}\frac{\left\|x_{k}-x_{k+1}\right\|^{2}}{2 \sigma}\nonumber\\
&+\left(1-{\sqrt{\tau\sigma}\|K\|}\right)\sum_{k=1}^{M-1}\frac{\left\|y_{k}-y_{k-1}\right\|^{2}}{2 \tau}+\frac{\left\|y_{M}-y_{M-1}\right\|^{2}}{2 \tau}\nonumber\\
&-\frac{\|y_{M}-y_{M-1}\|^2}{2\tau}-{\tau\sigma}\|K\|^2\frac{\|x_{M}-x\|^2}{2\sigma}\nonumber\\
&-\frac{L(T+1)^2}{2}\sum_{k=0}^{M-1}\|x_{k+1}-x_k\|^2\label{fs},
\end{align}
which can be written as
\begin{align}
\frac{\left\|x-x_{0}\right\|^2}{2 \sigma}+\frac{\left\|y-y_{0}\right\|^{2}}{2 \tau}
\geq&\sum_{k=1}^M\left(\mathcal{L}(x_{k},y)-\mathcal{L}(x,y_{k})\right)+\left(1-\tau\sigma\|K\|^2\right)\frac{\left\|x-x_{M}\right\|^{2}}{2 \sigma}\nonumber\\
&+\frac{\left\|y-y_{M}\right\|^{2}}{2 \tau}+\left(1-\sqrt{\tau\sigma}\|K\|-\sigma L(T+1)^2\right)\sum_{k=1}^{M}\frac{\left\|x_{k}-x_{k-1}\right\|^{2}}{2\sigma}\nonumber\\
&+\left(1-{\sqrt{\tau\sigma}\|K\|}\right)\sum_{k=1}^{M-1}\frac{\left\|y_{k}-y_{k-1}\right\|^{2}}{2 \tau} \label{f1}.
\end{align}
Taking $x=\hat x$ and $y=\hat y$, using \eqref{lieq} and condition \eqref{cond1}, we have
\begin{align}
\frac{\left\|\hat x-x_{0}\right\|}{2 \sigma}+\frac{\left\|\hat y-y_{0}\right\|^{2}}{2 \tau}\geq \left(1-\sigma\tau\|K\|^2\right)\left(\frac{\left\|\hat x-x_{M}\right\|^{2}}{2 \sigma}+\frac{\left\|\hat y-y_{M}\right\|^{2}}{2 \tau}\right),
\end{align}
which leads to statement (i).

For any $(x,y)\in\mathbb{R}^{d_1}\times \mathbb{R}^{d_2}$, it follows from \eqref{f1} that
\begin{align}
\frac{1}{M}\left(\frac{\left\|\hat x-x_{0}\right\|}{2 \sigma}+\frac{\left\|\hat y-y_{0}\right\|^{2}}{2 \tau}\right)\geq&\frac{1}{M}\left(\sum_{k=0}^{M-1}\mathcal{L}(x_{k+1},y)-\mathcal{L}(x,y_{k+1})\right)\nonumber\\
\geq&\mathcal{L}(\bar x_{M},y)-\mathcal{L}(x,\bar y_{M}),\label{f37}
\end{align}
where the convexity of $f$ and $h^\ast$ is employed in the second inequality.
Then, for any bounded closed sets $B_1$ and $B_2$, we have
\begin{align*}
\frac{1}{M}\max_{(x,y)\in B_1\times B_2}\left(\frac{\left\|x-x_{0}\right\|}{2 \sigma}+\frac{\left\|y-y_{0}\right\|^{2}}{2 \tau}\right)&\geq \max_{(x,y)\in B_1\times B_2}\left(\mathcal{L}(\bar x_M,y)-\mathcal{L}(x,\bar y_M)\right)\\
&=\mathcal{G}_{B_1\times B_2}(\bar x_M, \bar y_M).
\end{align*}
Suppose that $(x^\ast,y^\ast)$ is a cluster point of the sequence $\{(\bar x_k,\bar y_k)\}$. Since $f$ and $h^\ast$ are assumed to be closed, it follows from \eqref{f37} that
\begin{align*}
0\geq\mathcal{L}(x^\ast,y)-\mathcal{L}(x, y^\ast),~~~\forall (x,y)\in\mathbb{R}^{d_1}\times \mathbb{R}^{d_2},
\end{align*}
which implies $(x^\ast,y^\ast)$ is a saddle point according to the definition \eqref{equiv1}. Then statement (ii) is proved.

Since statement (i) implies $\{(x_k,y_k)\}$ is bounded, then there exists some subsequence $\{(x_{k_n},y_{k_n})\}$ converging to some point $(x^\ast,y^\ast)$. Taking $(x,y)=(\hat x,\hat y)$ in \eqref{f1}, we have $\lim_k(x_k-x_{k-1})=\lim_k(y_k-y_{k-\zeta})=0$ for any nonnegative integer $\zeta$. Then, $\{x_{k_n-\zeta}\}$ and $\{y_{k_n-1}\}$ converges to $x^\ast$ and $y^\ast$ respectively for any given $\zeta$, which follows that $(x^\ast,y^\ast)$ is a fixed point. Then $(x^\ast,y^\ast)$ is a saddle-point for $\mathcal{L}$. Taking $(x,y)=(x^\ast,y^\ast)$ in \eqref{fs491}, summing \eqref{fs491} from $k=k_n$ to $M-1$, we obtain 
\begin{align*}
\frac{\left\|x^\ast-x_{k_n}\right\|}{2 \sigma}+\frac{\left\|y^\ast-y_{k_n}\right\|^{2}}{2 \tau}
\geq&\frac{\left\|x^\ast-x_{M}\right\|^{2}}{2 \sigma}+\frac{\left\|y^\ast-y_{M}\right\|^{2}}{2 \tau}+(1-\sqrt{\tau\sigma}\|K\|)\sum_{k=k_n}^{M-1}\frac{\left\|y_{k}-y_{k+1}\right\|^{2}}{2 \sigma}\\
&+\left(1-{\sqrt{\tau\sigma}\|K\|}\right)\sum_{k=k_n}^{M-1}\frac{\left\|y_{k}-y_{k-1}\right\|^{2}}{2 \tau}+\frac{\left\|y_{M}-y_{M-1}\right\|^{2}}{2 \tau}-\frac{\left\|y_{k_n}-y_{k_n-1}\right\|^{2}}{2 \tau}\\
&+\left\langle K^T\left(y_{M}-y_{M-1}\right), x_{M}-x\right\rangle-\left\langle K^T\left(y_{k_n}-y_{k_n-1}\right), x_{k_n}-x_\ast\right\rangle\\
&-\frac{L(T+1)^2}{2}\sum_{k=k_n-T}^{M-1}\|x_{k+1}-x_k\|^2,
\end{align*}
which implies that the sequence $\{(x_M,y_M )\}$ converges to $(x^\ast,y^\ast)$. Statement (iii) is proved.
\qed
\subsection{Linear convergence}
In the second analysis, we choose the uncertain variable  $\bar{y} = y_k+\theta(y_k-y_{k-1})$. Under the strongly convexity of  $f$ and $h^\ast$, the linear convergence of iteration sequence can be obtained.
To obtain the linear convergence, we slightly modify a lemma from \cite{aytekin2016analysis}, where the nonnegativity of  $\{V_k\}$ is no longer required. The proof is omitted since the technical details are almost the same.
\begin{lemm}\label{lemm}Assume that the real sequence $\{V_k\}$ and the nonnegative sequence $\left\{\omega_{k}\right\}$ satisfy the following inequality:
\begin{align}
V_{k}\geq \frac{1}{a} V_{k+1}+b\omega_k-c\sum_{j=k-k_0}^k \omega_j
\end{align}
for some real numbers $a \in (0,1),b, c \geq 0$, and some positive integer $k_{0} $. Also, assume 
that $\omega_{k}=0$ for $k<0$ , and that the following holds:
\begin{align*}
\frac{c}{1-a} \frac{1-a^{k_{0}+1}}{a^{k_{0}}} \leq b.
\end{align*}
Then $V_{k} \leq a^{k} V_{0}$ for all $k \geq 0$.
\end{lemm}

\begin{theo}\label{linear}Assume A1-A3, B1 and B2 hold. Also assume $\mathcal{L}$ has a unique saddle point $(\hat x,\hat y)$. Take $\bar y = y_{k}+\theta(y_{k}-y_{k-1})$  in each iteration, where 
\begin{align*}
(\min \{3\delta\sigma/2,2\gamma\tau\}+1)^{-1}\leq\theta\leq 1. 
\end{align*}
For any sufficiently small $\sigma$ and $\tau$, the sequence $\{(x_k,y_k)\}$ is linearly convergent to $(\hat x,\hat y)$ with a rate of $O(\omega^{-k/2})$, where
\begin{equation*}
\omega:=\frac{1+\theta \sigma\|K\|}{\sigma\|K\|+1+\min \{3 \delta \sigma / 2,2 \gamma \tau\}}
\end{equation*}
\end{theo}
\textbf{Proof.}
Employing the $\delta_i$-strong convexity and $L_i$-smoothness of $F_i^\ast$, we have
\begin{align}
f_i (x)\geq& f_i({x_{k-\tau_k^i}})+\left\langle \nabla f_i(x_{k-\tau_k^i}), x-x_{k-\tau_k^i} \right\rangle+\frac{\delta_i}{2}\|x-x_{k-\tau_k^i}\|^2\nonumber\\
\geq& f_i({x_{k-\tau_k^i}})+\left\langle \nabla f_i(x_{k-\tau_k^i}), x-x_{k-\tau_k^i} \right\rangle+\frac{\delta_i}{4}\|x-x_{k+1}\|^2-\frac{\delta_i}{2}\|x_{k+1}-x_{k-\tau_k^i}\|^2\nonumber\\
\geq& f_i({x_{k+1}})+\left\langle \nabla f_i(x_{k-\tau_k^i}), x-x_{k+1} \right\rangle+\frac{\delta_i}{4}\|x-x_{k+1}\|^2-\frac{\delta_i+L_i}{2}\|x_{k+1}-x_{k-\tau_k^i}\|^2\nonumber\\
\geq& f_i({x_{k+1}})+\left\langle \nabla f_i(x_{k-\tau_k^i}), x-x_{k+1} \right\rangle+\frac{\delta_i}{4}\|x-x_{k+1}\|^2\nonumber\\
&-\frac{(L_i+\delta_i)(T+1)}{2}\sum_{j=k-T}^k\|x_{j+1}-x_{j}\|^2\label{f29},
\end{align}
where $2\|a\|^2\geq \|a+b\|^2-2\|b\|^2$ is used in the second inequality.
Summing \eqref{f29} from $i=1$ to $M$ yields
\begin{align}
f(x) \geq & f\left(x_{k+1}\right)+\left\langle g_k, x-x_{k+1}\right\rangle+\frac{\delta}{4}\|x-x_{k+1}\|^2\nonumber\\
&-\frac{(L+\delta)(T+1)}{2}\sum_{j=k-T}^k\|x_{j+1}-x_j\|^2 \nonumber\\
=&f\left(x_{k+1}\right)+\left\langle \frac{x_k-x_{k+1}}{\sigma}-K^T \overline{y}, x-x_{k+1}\right\rangle+\frac{\delta}{4}\|x-x_{k+1}\|^2\nonumber\\
&-\frac{(L+\delta)(T+1)}{2}\sum_{j=k-T}^k\|x_{j+1}-x_j\|^2.\label{fs1}
\end{align}
Recall the $\gamma$-strong convexity of $h^\ast$, we have
\begin{align} 
h^\ast(y) \geq h^\ast\left(y_{k+1}\right)&+\left\langle\frac{y_{k}-y_{k+1}}{\tau}, y-y_{k+1}\right\rangle\nonumber\\
&+\left\langle K^T\left(y-y_{k+1}\right), x_{k+1}\right\rangle+\frac{\gamma}{2}\left\|y-y_{k+1}\right\|^{2}.\label{gs2}
\end{align}
Since $\hat y$ is a minimizer of $\langle K \hat x, y\rangle+f(\hat x)-h^\ast(y)$ and $h^\ast$ is $\gamma$-strongly convex, we have
\begin{align}
-\langle K^T y, \hat{x}\rangle-f(\hat{x})+h^\ast(y)\geq -\langle K^T \hat y, \hat{x}\rangle-f(\hat{x})+h^\ast(\hat y)+\frac{\gamma}{2}\|y-\hat y\|^2.\label{tp1}
\end{align}
In the same way, we also obtain
\begin{align}
\langle K x, \hat y\rangle+ f(x)-h^\ast(\hat{y})\geq \langle K \hat x, \hat{y}\rangle+f(\hat{x})-h^\ast(\hat y)+\frac{\delta}{2}\|x-\hat x\|^2.\label{tp2}
\end{align}

The sum of \eqref{tp1} and \eqref{tp2} is
\begin{align}
\left[\langle K x, \hat y\rangle+f(x)-h^\ast(\hat y)\right]-\left[\langle K^T y, \hat x\rangle+f(\hat x)-h^\ast(y)\right]\geq \frac{\gamma}{2}\|y-\hat{y}\|^{2}+\frac{\delta}{2}\|x-\hat x\|^2.\label{cs3}
\end{align}
Combining \eqref{fs1}, \eqref{gs2} and \eqref{cs3} yields
\begin{align}
&\frac{\left\|\hat{x}-x_k\right\|^{2}}{2 \sigma}+\frac{\left\|\hat{y}-y_{k}\right\|^{2}}{2 \tau}\nonumber\\
\geq&\left( \frac{3\delta}{2}+\frac{1}{\sigma}\right)\frac{\left\|\hat{x}-x_{k+1}\right\|^{2}}{2}+\left(2 \gamma+\frac{1}{\tau}\right) \frac{\left\|\hat{y}-y_{k+1}\right\|^{2}}{2}\nonumber\\
&+\frac{\left\|x_k-x_{k+1}\right\|^{2}}{2 \sigma}+\frac{\left\|y_{k}-y_{k+1}\right\|^{2}}{2 \tau}-\frac{(L+\delta)(T+1)}{2}\sum_{j=k-T}^k\|x_{j+1}-x_{j}\|^2\nonumber\\
&-\left\langle K^T\left(y_{k+1}-\overline{y}\right), x_{k+1}-\hat{x}\right\rangle \label{crowd1}
\end{align}
Define $a:=\min\left\{1+3\delta\sigma/2,1+2\gamma\tau\right\}^{-1}$. Taking $\bar y=y_k+\theta(y_k-y_{k-1})$, and $$
\omega:=a \frac{1+\theta \sigma\|K\|}{1+a \sigma\|K\|}
$$, then we obtain:
\begin{align}
&-\left\langle K^T\left(y_{k+1}-\overline{y}\right), x_{k+1}-\hat{x}\right\rangle\\
=&-\left\langle K^T\left(y_{k+1}-y_{k}\right), x_{k+1}-\hat{x}\right\rangle+\theta\left\langle K^T\left(y_{k}-y_{k-1}\right), x_{k+1}-\hat{x}\right\rangle\nonumber\\
=&{-\left\langle K\left(y_{k+1}-y_{k}\right), x_{k+1}-\hat{x}\right\rangle+\omega\left\langle K^T\left(y_{k}-y_{k-1}\right), x_k-\hat{x}\right\rangle}\nonumber\\
&+\omega\left\langle K^T\left(y_{k}-y_{k-1}\right), x_{k+1}-x_k\right\rangle\nonumber \nonumber\\ 
&+(\theta-\omega)\left\langle K^T\left(y_{k}-y_{k-1}\right), x_{k+1}-\hat{x}\right\rangle\nonumber\\
\geq&-\left\langle K^T\left(y_{k+1}-y_{k}\right), x_{k+1}-\hat{x}\right\rangle+\omega\left\langle K^T\left(y_{k}-y_{k-1}\right), x_k-\hat{x}\right\rangle \nonumber\\
&-\omega \|K\|\left( \frac{\omega\left\|y_{k}-y_{k-1}\right\|^{2}}{2}+\frac{\left\|x_{k+1}-x_k\right\|^{2}}{2\omega }\right)\nonumber\\
&-(\theta-\omega) \|K\|\left(\frac{\omega\left\|y_{k}-y_{k-1}\right\|^{2}}{2}+\frac{\left\|x_{k+1}-\hat{x}\right\|^{2}}{2\omega}\right)\label{sub},
\end{align}
where the last inequality hold due to $\omega\leq\theta$.
Then by definition of $\omega$, 
$$
\frac{\omega}{\theta}=\frac{\frac{1}{\theta}+\sigma\|K\|}{\frac{1}{a}+\sigma\|K\|} \leq 1
$$
we get $a\leq\omega\leq\theta$.
Substituting \eqref{sub} into \eqref{crowd1} yields

\begin{align}
&\frac{\left\|\hat{x}-x_k\right\|^{2}}{2 \sigma}+\frac{\left\|\hat{y}-y_{k}\right\|^{2}}{2 \tau}\nonumber\\
\geq&\left( \frac{3\delta}{2}+\frac{1}{\sigma}\right)\frac{\left\|\hat{x}-x_{k+1}\right\|^{2}}{2}+\left(2 \gamma+\frac{1}{\tau}\right) \frac{\left\|\hat{y}-y_{k+1}\right\|^{2}}{2}\nonumber\\
&+\frac{\left\|x_k-x_{k+1}\right\|^{2}}{2 \sigma}+\frac{\left\|y_{k}-y_{k+1}\right\|^{2}}{2 \tau}-\frac{(L+\delta)(T+1)}{2}\sum_{j=k-T}^k\|x_{j+1}-x_{j}\|^2\nonumber\\
&-\left\langle K^T\left(y_{k+1}-y_{k}\right), x_{k+1}-\hat{x}\right\rangle+\omega\left\langle K^T\left(y_{k}-y_{k-1}\right), x_k-\hat{x}\right\rangle \nonumber\\
&-\theta\omega\|K\|\frac{\left\|y_{k}-y_{k-1}\right\|^{2}}{2}-\|K\|\frac{\|x_{k+1}-x_k\|^2}{2}-\frac{\theta-\omega}{\omega}\|K\|\frac{\|x_{k+1}-\hat x\|^2}{2}.\label{ref97}
\end{align}

Denote
\begin{align*}
V_k&:=\frac{\left\|\hat{x}-x_k\right\|^{2}}{2 \sigma}+\frac{\left\|\hat{y}-y_{k}\right\|^{2}}{2 \tau},\\
C_1&:=(L+\delta)(T+1),\\
\alpha_k &:=\frac{\|y_{k-1}-y_{k}\|^2}{2\tau},\beta_k :=\frac{\|x_k-x_{k+1}\|^2}{2\sigma},\\
\xi_k&:=-\langle K^T(y_k-y_{k-1}),x_k-\hat x\rangle,
\end{align*}
Then it follows from \eqref{ref97} that
\begin{align}
V_k\geq \frac{1}{\omega} V_{k+1}+(1-\|K\|\sigma)\beta_k-\sigma C_1\sum_{j=k-T}^k \beta_j+\xi_{k+1}-\omega \xi_{n}-\theta\omega \tau\|K\|\alpha_{k}+\alpha_{k+1}.\label{ccondc}
\end{align}
Note that for sufficiently  small $\tau$, we have $\theta\tau\|K\|\leq1$, which implies
\begin{align}
\left(V_k+\omega\xi_k+\omega\alpha_k\right)\geq \frac{1}{\omega} \left(V_{k+1}+\omega\xi_{k+1}+\omega\alpha_{k+1}\right)+(1-\|K\|\sigma) \beta_k-\sigma C_1\sum_{j=k-T}^k \beta_j\label{posb}
\end{align}
We verify the conditions of Lemma \ref{lemm}:
\begin{align}
\frac{\sigma C_1}{1-\omega}\cdot\frac{1-\omega^{T+1}}{\omega^T}&\leq \sigma C_1 (T+1) \omega^{-T}\leq\sigma C_1(T+1)a^{-T}\nonumber\\
&\leq \sigma C_1(T+1)\cdot\min\left\{1+3\delta \sigma/2,1+2\gamma\tau\right\}^T\nonumber\\
&\leq 1-\|K\|\sigma,\label{C1}
\end{align}
where the last inequality holds since we required sufficiently small $\sigma$ and $\tau$. It immediately follows from Lemma \ref{lemm} that
\begin{align}
V_k+\omega\xi_k+\omega \alpha_k\leq\omega^{k} V_0,~~~\forall k\geq 0.
\end{align}
Then we obtain

\begin{align}
\omega^{k}V_0
\geq& V_k-\omega \langle K^T(y_k-y_{k-1}),x_k-\hat x\rangle+\omega\frac{\|y_{k-1}-y_{k}\|^2}{2\tau}\nonumber\\
\geq & V_k-\omega\frac{\|y_{k-1}-y_k\|^2}{2\tau}-\omega\sigma\tau\|K\|^2\frac{\|x_k-\hat x\|^2}{2\sigma}+\omega\frac{\|y_{k-1}-y_{k}\|^2}{2\tau}\nonumber \\
\geq&\left(\frac{\|y_k-\hat y\|^2}{2\tau}+(1-\omega\sigma\tau\|K\|^2)\frac{\|x_k-\hat x\|^2}{2\sigma}\right)\nonumber\\
\geq&\left(\frac{\|y_k-\hat y\|^2}{2\tau}+(1-\sigma\tau\|K\|^2)\frac{\|x_k-\hat x\|^2}{2\sigma}\right).\label{linineq}
\end{align}
Therefore, the linear convergence of $\{(x_k,y_k)\}$ to the saddle point $(\hat x,\hat y)$ is proved.
\qed

In Theorem \ref{linear}, we required $\sigma$ and $\tau$ to be sufficiently small. In particular, we could choose the parameters satisfying that
\begin{align}
&\sigma(L+\delta)(T+1)\min\{1+3/2\cdot\gamma\tau,1+2\sigma\delta\}^T+\sigma\|K\|\leq 1,\label{conda}\\
&\sigma\tau\|K\|^2<1,\label{condb}\\
&\theta\tau\|K\|\leq1, \label{condc}
\end{align}
where \eqref{conda} is required by \eqref{C1}, \eqref{condb} is required by \eqref{linineq}, and \eqref{condc} is required by \eqref{ccondc}. Also note that \eqref{conda} implies $\sigma\|K\|< 1$, which is used in \eqref{posb} to guarantee the positive of second part in the last formula.
\subsection{ The Sublinear Convergence in AH formation}

The Arrow-Hurwicz (AH) method has been studied in \cite{nedic2009subgradient}, which has a relatively fast rate of convergence in Primal-Dual problem. Like AH formation, take  $\bar y=y_{k}$ in PD-PIAG.
\begin{theo} Assume that A1-A3, B1 hold, the problem has a saddle point $(\hat x,\hat y)$, and $D_x:=\sup_k\|x_k-\hat x\|<\infty$. Take $\bar y = y_{k}$ in each iteration. Choose $\sigma$ and $\tau$ such that  $1-\sigma L(T +1)^2<0,\|K\|^2\tau\leq \delta$.
Then, 
\begin{enumerate}
\item [\rm(i)]
the sequence $\{(x_k,y_k)\}$ is bounded since
\begin{align}
\frac{\left\|x_{k}-\hat{x}\right\|^{2}}{2 \sigma}+\frac{\left\|y_{k}-\hat{y}\right\|^{2}}{2 \tau}
\quad \leq C\left(\frac{\left\|x_{0}-\hat{x}\right\|^{2}}{2 \sigma}+\frac{\left\|y_{0}-\hat{y}\right\|^{2}}{2 \tau}\right),
\end{align}
where the constant $C=(1-\tau\sigma\|K\|^2)^{-1}$. 
\item[\rm(ii)]
Define the averaged sequences $\bar x_M=(\sum_{k=1}^M x_k)/M$ and $\bar y_M=(\sum_{k=1}^M y_k)/M$ for all $M>0$.
Then for any bounded closed set $B_1\times B_2\subset X\times Y$, the restricted gap has the following bound
\begin{align}
\mathcal{G}_{B_{1} \times B_{2}}\left(\bar x_{M}, \bar y_{M}\right) \leq \frac{1}{M}\max _{(x, y) \in B_{1} \times B_{2}} \left\{\frac{\left\|x-x_{0}\right\|^{2}}{2 \sigma}+\frac{\left\|y-y_{0}\right\|^{2}}{2 \tau}\right\}.
\end{align}
Moreover, the cluster points of $\{(x_M,y_M)\}$ are saddle points. 
\item[\rm(iii)]
There exists a saddle point $(x^\ast,y^\ast)$ such that $x_k\rightarrow x^\ast$ and $y_k\rightarrow y^\ast$.
\end{enumerate}
\end{theo}
\textbf{Proof.}
Taking $ \bar y=y_{k}$ in \eqref{algorithm1}.  Similar to the proof of \eqref{fs}, we have
\begin{align}
\frac{\left\|x-x_{k}\right\|^2}{2 \sigma}+\frac{\left\|y-y_{k}\right\|^{2}}{2 \tau}
\geq&\mathcal{L}(x_{k+1},y)-\mathcal{L}(x,y_{k+1})\nonumber+\left( \delta+\frac{1}{\sigma}\right)\frac{\left\|x-x_{k+1}\right\|^{2}}{2}\\
&+\frac{\left\|y-y_{k+1}\right\|^{2}}{2 \tau}+\frac{\left\|x_{k}-x_{k+1}\right\|^{2}}{2 \sigma}+\frac{\left\|y_{k}-y_{k+1}\right\|^{2}}{2 \tau}\nonumber\\
&+\left\langle K\left(x_{k+1}-x\right), y_{k+1}-y_k\right\rangle\nonumber\\
&-\sigma L(T+1)\sum_{j=k-T}^k\frac{\|x_{j+1}-x_j\|^2}{2\sigma}.\label{fs78}
\end{align}
Employing
\begin{align*}
\left\langle K\left(x_{k+1}-x\right), y_{k+1}- y_k\right\rangle&\geq-\|K\|^2\frac{\left\|y_{k+1}-y_{k}\right\|^{2}}{2 \delta}-\delta \frac{\left\|x_{k+1}- x\right\|^{2}}{2 },
\end{align*}
summing \eqref{fs78} from $K=0$ to $M-1$, we obtain
\begin{align}
&\sum_{k=1}^{M}\left[\mathcal{L}(x_k, y)-\mathcal{L}(x,y_k)\right] +\frac{\left\| x-x_{M}\right\|^{2}}{2 \sigma}+\frac{\left\| y-y_{M}\right\|^{2}}{2 \tau}\nonumber\\
&+(1-\sigma L(T+1)^2)\sum_{k=1}^{M} \frac{\left\|x_k-x_{k-1}\right\|^{2}}{2 \sigma}+\left(1-\frac{\|K\|^2\tau}{\delta}\right)\sum_{k=1}^{M} \frac{\left\|y_{k}-y_{k-1}\right\|^{2}}{2 \tau}\nonumber\\ 
&\leq  \frac{\left\|x-x_{0}\right\|^{2}}{2 \sigma}+\frac{\left\|y-y_{0}\right\|^{2}}{2 \tau}.\label{fs49}
\end{align}
Note that $\mathcal{L}(x_k, \hat y)-\mathcal{L}( \hat x,y_k)\geq 0$ due to \eqref{lieq}. Then Statement (i) is proved. 

It immediately follows from \eqref{fs49} that
\begin{align}
\mathcal{G}_{B_1\times B_2}(\bar x_k,\bar y_k)=&\max_{y\in B_2}\mathcal{L}(\bar x_k,y)-\min_{x\in B_1}\mathcal{L}(x,\bar y_k)\nonumber\\
\leq& \frac{1}{M}\sup_{(x,y)\in B_1\times B_2}\left\{\frac{\left\|x-x_{0}\right\|^{2}}{2 \sigma}+\frac{\left\|y-y_{0}\right\|^{2}}{2 \tau}\right\}.
\end{align}
holds for any closed bounded $B_1$ and $B_2$. Analogous  to the proof of Statement (ii) in Theorem \ref{theo1}, we prove any cluster point of $\{(\bar x_k,y_k)\}$ is a saddle point. Again, similar to the proof of Statement (iii) in Theorem \ref{theo1}, we prove $\{(x_k,y_k)\}$ converges to some saddle point $(x^\ast,y^\ast)$.
\qed


\section{Conclusion}
In this paper, we have analyzed the convergence of PD-PIAG for saddle-point problem. First, we proposed the PD-PIAG with extrapolations, and provided the sublinear convergence of its partial primal-dual gap, as well as the convergence of iterates to some saddle point. With strongly convexity of $f$ and $h^\ast$, we proved that the generated sequence $\{(x_k,y_k)\}$ converges to the saddle point. Then, we proposed PD-PIAG, the primal-dual gap  of which is proved to be sublinearly convergent under strong convexity of $f$. The proposed incremental aggregated methods can be reviewed as asynchronous variants of several existing primal-dual methods. However, a generalized framework should be established to analyze PD methods with incremental aggregated settings, which deserves further study.
\bibliographystyle{plain}
\bibliography{iag_refs}
\end{document}